# SIMPLE DIVISIBILITY RULES FOR THE 1$^{st}$ 1000 PRIME NUMBERS

C. C. Briggs
*Center for Academic Computing, Penn State University, University Park, PA 16802*
Saturday, December 25, 1999

**Abstract.** Simple divisibility rules are given for the 1$^{st}$ 1000 prime numbers.

This paper presents simple divisibility rules for the 1$^{st}$ 1000 prime numbers.

With two exceptions, the rules in question are based on the observation that if $M$ is an integer not divisible by 2 or by 5, then $M$ divides an arbitrary integer $N$ if it divides

$$N' = \frac{N - N \bmod 10}{10} + (m + A \times M) \times (N \bmod 10) \qquad (1)$$

$$= \frac{N - M \times (m' - 10 \times A) \times (N \bmod 10)}{10},$$

where

$$m = \frac{1 - m' \times M}{10}, \qquad (2)$$

$A$ is any integer (positive, negative, or 0), and, for $M > 0$,

$$m' = [3 \times (M \bmod 10) - 2 \times M] \bmod 16, \qquad (3)$$

the integer $A$ in Eq. (1) making it possible to customize the rule.

The rules for divisibility by 2 and by 5—the two aforesaid exceptions—are based on the observation that, for $p$ an arbitrary non-negative integer, $2^p$ (and, similarly, $5^p$) divides $N$ if it divides $N \bmod 10^p$, since $10^p$ divides $N \bmod 10^{p+1} - N \bmod 10^p$, while $2^p$ (and, similarly, $5^p$) divides $10^p = (2 \times 5)^p = 2^p \times 5^p$.

Simple divisibility rules for the 1$^{st}$ 1000 prime numbers based on the preceding observations (with $A = 0$ and $p = 1$) appear in Table 1 (see below), where—technically speaking—"number" means "integer," "times" means "multiplied by," and, for $N$ being "the number," "the last digit of the number" is given by $N \bmod 10$ and "the rest of the number" is given by $(N - N \bmod 10)/10$.

TABLE 1. SIMPLE DIVISIBILITY RULES FOR THE 1$^{st}$ 1000 PRIME NUMBERS

| # | | $M$ | | coef. | | | |
|---|---|---|---|---|---|---|---|
| 1. | A number is divisible by | 2 | if the last digit of the number is divisible by 2. | | | | |
| 2. | A number is divisible by | 3 | if | − 2 | times the last digit of the number added to the rest of the number is divisible by | 3. |
| 3. | A number is divisible by | 5 | if the last digit of the number is divisible by 5. | | | | |
| 4. | A number is divisible by | 7 | if | − 2 | times the last digit of the number added to the rest of the number is divisible by | 7. |
| 5. | A number is divisible by | 11 | if | − 1 | times the last digit of the number added to the rest of the number is divisible by | 11. |
| 6. | A number is divisible by | 13 | if | − 9 | times the last digit of the number added to the rest of the number is divisible by | 13. |
| 7. | " | 17 | " | − 5 | " | 17. |
| 8. | " | 19 | " | − 17 | " | 19. |
| 9. | " | 23 | " | − 16 | " | 23. |
| 10. | " | 29 | " | − 26 | " | 29. |
| 11. | A number is divisible by | 31 | if | − 3 | times the last digit of the number added to the rest of the number is divisible by | 31. |
| 12. | " | 37 | " | − 11 | " | 37. |
| 13. | " | 41 | " | − 4 | " | 41. |
| 14. | " | 43 | " | − 30 | " | 43. |
| 15. | " | 47 | " | − 14 | " | 47. |
| 16. | A number is divisible by | 53 | if | − 37 | times the last digit of the number added to the rest of the number is divisible by | 53. |
| 17. | " | 59 | " | − 53 | " | 59. |
| 18. | " | 61 | " | − 6 | " | 61. |
| 19. | " | 67 | " | − 20 | " | 67. |
| 20. | " | 71 | " | − 7 | " | 71. |
| 21. | A number is divisible by | 73 | if | − 51 | times the last digit of the number added to the rest of the number is divisible by | 73. |
| 22. | " | 79 | " | − 71 | " | 79. |
| 23. | " | 83 | " | − 58 | " | 83. |
| 24. | " | 89 | " | − 80 | " | 89. |
| 25. | " | 97 | " | − 29 | " | 97. |
| 26. | A number is divisible by | 101 | if | − 10 | times the last digit of the number added to the rest of the number is divisible by | 101. |
| 27. | " | 103 | " | − 72 | " | 103. |
| 28. | " | 107 | " | − 32 | " | 107. |
| 29. | " | 109 | " | − 98 | " | 109. |
| 30. | " | 113 | " | − 79 | " | 113. |
| 31. | A number is divisible by | 127 | if | − 38 | times the last digit of the number added to the rest of the number is divisible by | 127. |
| 32. | " | 131 | " | − 13 | " | 131. |
| 33. | " | 137 | " | − 41 | " | 137. |
| 34. | " | 139 | " | − 125 | " | 139. |
| 35. | " | 149 | " | − 134 | " | 149. |
| 36. | A number is divisible by | 151 | if | − 15 | times the last digit of the number added to the rest of the number is divisible by | 151. |
| 37. | " | 157 | " | − 47 | " | 157. |
| 38. | " | 163 | " | − 114 | " | 163. |
| 39. | " | 167 | " | − 50 | " | 167. |
| 40. | " | 173 | " | − 121 | " | 173. |
| 41. | A number is divisible by | 179 | if | − 161 | times the last digit of the number added to the rest of the number is divisible by | 179. |
| 42. | " | 181 | " | − 18 | " | 181. |
| 43. | " | 191 | " | − 19 | " | 191. |
| 44. | " | 193 | " | − 135 | " | 193. |
| 45. | " | 197 | " | − 59 | " | 197. |

TABLE 1. SIMPLE DIVISIBILITY RULES FOR THE 1$^{st}$ 1000 PRIME NUMBERS—*Continued*

| | | | | | |
|---|---|---|---|---|---|
| 46. | A number is divisible by | 199 if | − 179 times the last digit of the number added to the rest of the number is divisible by | 199. |
| 47. | " | 211 " | − 21 | " | 211. |
| 48. | " | 223 " | − 156 | " | 223. |
| 49. | " | 227 " | − 68 | " | 227. |
| 50. | " | 229 " | − 206 | " | 229. |
| 51. | A number is divisible by | 233 if | − 163 times the last digit of the number added to the rest of the number is divisible by | 233. |
| 52. | " | 239 " | − 215 | " | 239. |
| 53. | " | 241 " | − 24 | " | 241. |
| 54. | " | 251 " | − 25 | " | 251. |
| 55. | " | 257 " | − 77 | " | 257. |
| 56. | A number is divisible by | 263 if | − 184 times the last digit of the number added to the rest of the number is divisible by | 263. |
| 57. | " | 269 " | − 242 | " | 269. |
| 58. | " | 271 " | − 27 | " | 271. |
| 59. | " | 277 " | − 83 | " | 277. |
| 60. | " | 281 " | − 28 | " | 281. |
| 61. | A number is divisible by | 283 if | − 198 times the last digit of the number added to the rest of the number is divisible by | 283. |
| 62. | " | 293 " | − 205 | " | 293. |
| 63. | " | 307 " | − 92 | " | 307. |
| 64. | " | 311 " | − 31 | " | 311. |
| 65. | " | 313 " | − 219 | " | 313. |
| 66. | A number is divisible by | 317 if | − 95 times the last digit of the number added to the rest of the number is divisible by | 317. |
| 67. | " | 331 " | − 33 | " | 331. |
| 68. | " | 337 " | − 101 | " | 337. |
| 69. | " | 347 " | − 104 | " | 347. |
| 70. | " | 349 " | − 314 | " | 349. |
| 71. | A number is divisible by | 353 if | − 247 times the last digit of the number added to the rest of the number is divisible by | 353. |
| 72. | " | 359 " | − 323 | " | 359. |
| 73. | " | 367 " | − 110 | " | 367. |
| 74. | " | 373 " | − 261 | " | 373. |
| 75. | " | 379 " | − 341 | " | 379. |
| 76. | A number is divisible by | 383 if | − 268 times the last digit of the number added to the rest of the number is divisible by | 383. |
| 77. | " | 389 " | − 350 | " | 389. |
| 78. | " | 397 " | − 119 | " | 397. |
| 79. | " | 401 " | − 40 | " | 401. |
| 80. | " | 409 " | − 368 | " | 409. |
| 81. | A number is divisible by | 419 if | − 377 times the last digit of the number added to the rest of the number is divisible by | 419. |
| 82. | " | 421 " | − 42 | " | 421. |
| 83. | " | 431 " | − 43 | " | 431. |
| 84. | " | 433 " | − 303 | " | 433. |
| 85. | " | 439 " | − 395 | " | 439. |
| 86. | A number is divisible by | 443 if | − 310 times the last digit of the number added to the rest of the number is divisible by | 443. |
| 87. | " | 449 " | − 404 | " | 449. |
| 88. | " | 457 " | − 137 | " | 457. |
| 89. | " | 461 " | − 46 | " | 461. |
| 90. | " | 463 " | − 324 | " | 463. |
| 91. | A number is divisible by | 467 if | − 140 times the last digit of the number added to the rest of the number is divisible by | 467. |
| 92. | " | 479 " | − 431 | " | 479. |
| 93. | " | 487 " | − 146 | " | 487. |
| 94. | " | 491 " | − 49 | " | 491. |
| 95. | " | 499 " | − 449 | " | 499. |
| 96. | A number is divisible by | 503 if | − 352 times the last digit of the number added to the rest of the number is divisible by | 503. |
| 97. | " | 509 " | − 458 | " | 509. |
| 98. | " | 521 " | − 52 | " | 521. |
| 99. | " | 523 " | − 366 | " | 523. |
| 100. | " | 541 " | − 54 | " | 541. |
| 101. | A number is divisible by | 547 if | − 164 times the last digit of the number added to the rest of the number is divisible by | 547. |
| 102. | " | 557 " | − 167 | " | 557. |
| 103. | " | 563 " | − 394 | " | 563. |
| 104. | " | 569 " | − 512 | " | 569. |
| 105. | " | 571 " | − 57 | " | 571. |
| 106. | A number is divisible by | 577 if | − 173 times the last digit of the number added to the rest of the number is divisible by | 577. |
| 107. | " | 587 " | − 176 | " | 587. |
| 108. | " | 593 " | − 415 | " | 593. |
| 109. | " | 599 " | − 539 | " | 599. |
| 110. | " | 601 " | − 60 | " | 601. |
| 111. | A number is divisible by | 607 if | − 182 times the last digit of the number added to the rest of the number is divisible by | 607. |
| 112. | " | 613 " | − 429 | " | 613. |
| 113. | " | 617 " | − 185 | " | 617. |
| 114. | " | 619 " | − 557 | " | 619. |
| 115. | " | 631 " | − 63 | " | 631. |
| 116. | A number is divisible by | 641 if | − 64 times the last digit of the number added to the rest of the number is divisible by | 641. |
| 117. | " | 643 " | − 450 | " | 643. |
| 118. | " | 647 " | − 194 | " | 647. |
| 119. | " | 653 " | − 457 | " | 653. |
| 120. | " | 659 " | − 593 | " | 659. |

TABLE 1. SIMPLE DIVISIBILITY RULES FOR THE $1^{st}$ 1000 PRIME NUMBERS—*Continued*

| | | | | | |
|---|---|---|---|---|---|
| 121. | A number is divisible by | 661 if | $-66$ times the last digit of the number added to the rest of the number is divisible by | 661. |
| 122. | " | 673 " | $-471$ | " | 673. |
| 123. | " | 677 " | $-203$ | " | 677. |
| 124. | " | 683 " | $-478$ | " | 683. |
| 125. | " | 691 " | $-69$ | " | 691. |
| 126. | A number is divisible by | 701 if | $-70$ times the last digit of the number added to the rest of the number is divisible by | 701. |
| 127. | " | 709 " | $-638$ | " | 709. |
| 128. | " | 719 " | $-647$ | " | 719. |
| 129. | " | 727 " | $-218$ | " | 727. |
| 130. | " | 733 " | $-513$ | " | 733. |
| 131. | A number is divisible by | 739 if | $-665$ times the last digit of the number added to the rest of the number is divisible by | 739. |
| 132. | " | 743 " | $-520$ | " | 743. |
| 133. | " | 751 " | $-75$ | " | 751. |
| 134. | " | 757 " | $-227$ | " | 757. |
| 135. | " | 761 " | $-76$ | " | 761. |
| 136. | A number is divisible by | 769 if | $-692$ times the last digit of the number added to the rest of the number is divisible by | 769. |
| 137. | " | 773 " | $-541$ | " | 773. |
| 138. | " | 787 " | $-236$ | " | 787. |
| 139. | " | 797 " | $-239$ | " | 797. |
| 140. | " | 809 " | $-728$ | " | 809. |
| 141. | A number is divisible by | 811 if | $-81$ times the last digit of the number added to the rest of the number is divisible by | 811. |
| 142. | " | 821 " | $-82$ | " | 821. |
| 143. | " | 823 " | $-576$ | " | 823. |
| 144. | " | 827 " | $-248$ | " | 827. |
| 145. | " | 829 " | $-746$ | " | 829. |
| 146. | A number is divisible by | 839 if | $-755$ times the last digit of the number added to the rest of the number is divisible by | 839. |
| 147. | " | 853 " | $-597$ | " | 853. |
| 148. | " | 857 " | $-257$ | " | 857. |
| 149. | " | 859 " | $-773$ | " | 859. |
| 150. | " | 863 " | $-604$ | " | 863. |
| 151. | A number is divisible by | 877 if | $-263$ times the last digit of the number added to the rest of the number is divisible by | 877. |
| 152. | " | 881 " | $-88$ | " | 881. |
| 153. | " | 883 " | $-618$ | " | 883. |
| 154. | " | 887 " | $-266$ | " | 887. |
| 155. | " | 907 " | $-272$ | " | 907. |
| 156. | A number is divisible by | 911 if | $-91$ times the last digit of the number added to the rest of the number is divisible by | 911. |
| 157. | " | 919 " | $-827$ | " | 919. |
| 158. | " | 929 " | $-836$ | " | 929. |
| 159. | " | 937 " | $-281$ | " | 937. |
| 160. | " | 941 " | $-94$ | " | 941. |
| 161. | A number is divisible by | 947 if | $-284$ times the last digit of the number added to the rest of the number is divisible by | 947. |
| 162. | " | 953 " | $-667$ | " | 953. |
| 163. | " | 967 " | $-290$ | " | 967. |
| 164. | " | 971 " | $-97$ | " | 971. |
| 165. | " | 977 " | $-293$ | " | 977. |
| 166. | A number is divisible by | 983 if | $-688$ times the last digit of the number added to the rest of the number is divisible by | 983. |
| 167. | " | 991 " | $-99$ | " | 991. |
| 168. | " | 997 " | $-299$ | " | 997. |
| 169. | " | 1009 " | $-908$ | " | 1009. |
| 170. | " | 1013 " | $-709$ | " | 1013. |
| 171. | A number is divisible by | 1019 if | $-917$ times the last digit of the number added to the rest of the number is divisible by | 1019. |
| 172. | " | 1021 " | $-102$ | " | 1021. |
| 173. | " | 1031 " | $-103$ | " | 1031. |
| 174. | " | 1033 " | $-723$ | " | 1033. |
| 175. | " | 1039 " | $-935$ | " | 1039. |
| 176. | A number is divisible by | 1049 if | $-944$ times the last digit of the number added to the rest of the number is divisible by | 1049. |
| 177. | " | 1051 " | $-105$ | " | 1051. |
| 178. | " | 1061 " | $-106$ | " | 1061. |
| 179. | " | 1063 " | $-744$ | " | 1063. |
| 180. | " | 1069 " | $-962$ | " | 1069. |
| 181. | A number is divisible by | 1087 if | $-326$ times the last digit of the number added to the rest of the number is divisible by | 1087. |
| 182. | " | 1091 " | $-109$ | " | 1091. |
| 183. | " | 1093 " | $-765$ | " | 1093. |
| 184. | " | 1097 " | $-329$ | " | 1097. |
| 185. | " | 1103 " | $-772$ | " | 1103. |
| 186. | A number is divisible by | 1109 if | $-998$ times the last digit of the number added to the rest of the number is divisible by | 1109. |
| 187. | " | 1117 " | $-335$ | " | 1117. |
| 188. | " | 1123 " | $-786$ | " | 1123. |
| 189. | " | 1129 " | $-1016$ | " | 1129. |
| 190. | " | 1151 " | $-115$ | " | 1151. |
| 191. | A number is divisible by | 1153 if | $-807$ times the last digit of the number added to the rest of the number is divisible by | 1153. |
| 192. | " | 1163 " | $-814$ | " | 1163. |
| 193. | " | 1171 " | $-117$ | " | 1171. |
| 194. | " | 1181 " | $-118$ | " | 1181. |
| 195. | " | 1187 " | $-356$ | " | 1187. |

TABLE 1. SIMPLE DIVISIBILITY RULES FOR THE $1^{st}$ 1000 PRIME NUMBERS—*Continued*

196. A number is divisible by 1193 if − 835 times the last digit of the number added to the rest of the number is divisible by 1193.
197. " 1201 " − 120 " 1201.
198. " 1213 " − 849 " 1213.
199. " 1217 " − 365 " 1217.
200. " 1223 " − 856 " 1223.
201. A number is divisible by 1229 if − 1106 times the last digit of the number added to the rest of the number is divisible by 1229.
202. " 1231 " − 123 " 1231.
203. " 1237 " − 371 " 1237.
204. " 1249 " − 1124 " 1249.
205. " 1259 " − 1133 " 1259.
206. A number is divisible by 1277 if − 383 times the last digit of the number added to the rest of the number is divisible by 1277.
207. " 1279 " − 1151 " 1279.
208. " 1283 " − 898 " 1283.
209. " 1289 " − 1160 " 1289.
210. " 1291 " − 129 " 1291.
211. A number is divisible by 1297 if − 389 times the last digit of the number added to the rest of the number is divisible by 1297.
212. " 1301 " − 130 " 1301.
213. " 1303 " − 912 " 1303.
214. " 1307 " − 392 " 1307.
215. " 1319 " − 1187 " 1319.
216. A number is divisible by 1321 if − 132 times the last digit of the number added to the rest of the number is divisible by 1321.
217. " 1327 " − 398 " 1327.
218. " 1361 " − 136 " 1361.
219. " 1367 " − 410 " 1367.
220. " 1373 " − 961 " 1373.
221. A number is divisible by 1381 if − 138 times the last digit of the number added to the rest of the number is divisible by 1381.
222. " 1399 " − 1259 " 1399.
223. " 1409 " − 1268 " 1409.
224. " 1423 " − 996 " 1423.
225. " 1427 " − 428 " 1427.
226. A number is divisible by 1429 if − 1286 times the last digit of the number added to the rest of the number is divisible by 1429.
227. " 1433 " − 1003 " 1433.
228. " 1439 " − 1295 " 1439.
229. " 1447 " − 434 " 1447.
230. " 1451 " − 145 " 1451.
231. A number is divisible by 1453 if − 1017 times the last digit of the number added to the rest of the number is divisible by 1453.
232. " 1459 " − 1313 " 1459.
233. " 1471 " − 147 " 1471.
234. " 1481 " − 148 " 1481.
235. " 1483 " − 1038 " 1483.
236. A number is divisible by 1487 if − 446 times the last digit of the number added to the rest of the number is divisible by 1487.
237. " 1489 " − 1340 " 1489.
238. " 1493 " − 1045 " 1493.
239. " 1499 " − 1349 " 1499.
240. " 1511 " − 151 " 1511.
241. A number is divisible by 1523 if − 1066 times the last digit of the number added to the rest of the number is divisible by 1523.
242. " 1531 " − 153 " 1531.
243. " 1543 " − 1080 " 1543.
244. " 1549 " − 1394 " 1549.
245. " 1553 " − 1087 " 1553.
246. A number is divisible by 1559 if − 1403 times the last digit of the number added to the rest of the number is divisible by 1559.
247. " 1567 " − 470 " 1567.
248. " 1571 " − 157 " 1571.
249. " 1579 " − 1421 " 1579.
250. " 1583 " − 1108 " 1583.
251. A number is divisible by 1597 if − 479 times the last digit of the number added to the rest of the number is divisible by 1597.
252. " 1601 " − 160 " 1601.
253. " 1607 " − 482 " 1607.
254. " 1609 " − 1448 " 1609.
255. " 1613 " − 1129 " 1613.
256. A number is divisible by 1619 if − 1457 times the last digit of the number added to the rest of the number is divisible by 1619.
257. " 1621 " − 162 " 1621.
258. " 1627 " − 488 " 1627.
259. " 1637 " − 491 " 1637.
260. " 1657 " − 497 " 1657.
261. A number is divisible by 1663 if − 1164 times the last digit of the number added to the rest of the number is divisible by 1663.
262. " 1667 " − 500 " 1667.
263. " 1669 " − 1502 " 1669.
264. " 1693 " − 1185 " 1693.
265. " 1697 " − 509 " 1697.
266. A number is divisible by 1699 if − 1529 times the last digit of the number added to the rest of the number is divisible by 1699.
267. " 1709 " − 1538 " 1709.
268. " 1721 " − 172 " 1721.
269. " 1723 " − 1206 " 1723.
270. " 1733 " − 1213 " 1733.

TABLE 1. SIMPLE DIVISIBILITY RULES FOR THE 1$^{st}$ 1000 PRIME NUMBERS—*Continued*

271. A number is divisible by 1741 if  − 174 times the last digit of the number added to the rest of the number is divisible by 1741.
272. " 1747 " − 524 " 1747.
273. " 1753 " − 1227 " 1753.
274. " 1759 " − 1583 " 1759.
275. " 1777 " − 533 " 1777.
276. A number is divisible by 1783 if  − 1248 times the last digit of the number added to the rest of the number is divisible by 1783.
277. " 1787 " − 536 " 1787.
278. " 1789 " − 1610 " 1789.
279. " 1801 " − 180 " 1801.
280. " 1811 " − 181 " 1811.
281. A number is divisible by 1823 if  − 1276 times the last digit of the number added to the rest of the number is divisible by 1823.
282. " 1831 " − 183 " 1831.
283. " 1847 " − 554 " 1847.
284. " 1861 " − 186 " 1861.
285. " 1867 " − 560 " 1867.
286. A number is divisible by 1871 if  − 187 times the last digit of the number added to the rest of the number is divisible by 1871.
287. " 1873 " − 1311 " 1873.
288. " 1877 " − 563 " 1877.
289. " 1879 " − 1691 " 1879.
290. " 1889 " − 1700 " 1889.
291. A number is divisible by 1901 if  − 190 times the last digit of the number added to the rest of the number is divisible by 1901.
292. " 1907 " − 572 " 1907.
293. " 1913 " − 1339 " 1913.
294. " 1931 " − 193 " 1931.
295. " 1933 " − 1353 " 1933.
296. A number is divisible by 1949 if  − 1754 times the last digit of the number added to the rest of the number is divisible by 1949.
297. " 1951 " − 195 " 1951.
298. " 1973 " − 1381 " 1973.
299. " 1979 " − 1781 " 1979.
300. " 1987 " − 596 " 1987.
301. A number is divisible by 1993 if  − 1395 times the last digit of the number added to the rest of the number is divisible by 1993.
302. " 1997 " − 599 " 1997.
303. " 1999 " − 1799 " 1999.
304. " 2003 " − 1402 " 2003.
305. " 2011 " − 201 " 2011.
306. A number is divisible by 2017 if  − 605 times the last digit of the number added to the rest of the number is divisible by 2017.
307. " 2027 " − 608 " 2027.
308. " 2029 " − 1826 " 2029.
309. " 2039 " − 1835 " 2039.
310. " 2053 " − 1437 " 2053.
311. A number is divisible by 2063 if  − 1444 times the last digit of the number added to the rest of the number is divisible by 2063.
312. " 2069 " − 1862 " 2069.
313. " 2081 " − 208 " 2081.
314. " 2083 " − 1458 " 2083.
315. " 2087 " − 626 " 2087.
316. A number is divisible by 2089 if  − 1880 times the last digit of the number added to the rest of the number is divisible by 2089.
317. " 2099 " − 1889 " 2099.
318. " 2111 " − 211 " 2111.
319. " 2113 " − 1479 " 2113.
320. " 2129 " − 1916 " 2129.
321. A number is divisible by 2131 if  − 213 times the last digit of the number added to the rest of the number is divisible by 2131.
322. " 2137 " − 641 " 2137.
323. " 2141 " − 214 " 2141.
324. " 2143 " − 1500 " 2143.
325. " 2153 " − 1507 " 2153.
326. A number is divisible by 2161 if  − 216 times the last digit of the number added to the rest of the number is divisible by 2161.
327. " 2179 " − 1961 " 2179.
328. " 2203 " − 1542 " 2203.
329. " 2207 " − 662 " 2207.
330. " 2213 " − 1549 " 2213.
331. A number is divisible by 2221 if  − 222 times the last digit of the number added to the rest of the number is divisible by 2221.
332. " 2237 " − 671 " 2237.
333. " 2239 " − 2015 " 2239.
334. " 2243 " − 1570 " 2243.
335. " 2251 " − 225 " 2251.
336. A number is divisible by 2267 if  − 680 times the last digit of the number added to the rest of the number is divisible by 2267.
337. " 2269 " − 2042 " 2269.
338. " 2273 " − 1591 " 2273.
339. " 2281 " − 228 " 2281.
340. " 2287 " − 686 " 2287.
341. A number is divisible by 2293 if  − 1605 times the last digit of the number added to the rest of the number is divisible by 2293.
342. " 2297 " − 689 " 2297.
343. " 2309 " − 2078 " 2309.
344. " 2311 " − 231 " 2311.
345. " 2333 " − 1633 " 2333.

TABLE 1. SIMPLE DIVISIBILITY RULES FOR THE 1$^{st}$ 1000 PRIME NUMBERS—*Continued*

346. A number is divisible by 2339 if − 2105 times the last digit of the number added to the rest of the number is divisible by 2339.
347. " 2341 " − 234 " 2341.
348. " 2347 " − 704 " 2347.
349. " 2351 " − 235 " 2351.
350. " 2357 " − 707 " 2357.
351. A number is divisible by 2371 if − 237 times the last digit of the number added to the rest of the number is divisible by 2371.
352. " 2377 " − 713 " 2377.
353. " 2381 " − 238 " 2381.
354. " 2383 " − 1668 " 2383.
355. " 2389 " − 2150 " 2389.
356. A number is divisible by 2393 if − 1675 times the last digit of the number added to the rest of the number is divisible by 2393.
357. " 2399 " − 2159 " 2399.
358. " 2411 " − 241 " 2411.
359. " 2417 " − 725 " 2417.
360. " 2423 " − 1696 " 2423.
361. A number is divisible by 2437 if − 731 times the last digit of the number added to the rest of the number is divisible by 2437.
362. " 2441 " − 244 " 2441.
363. " 2447 " − 734 " 2447.
364. " 2459 " − 2213 " 2459.
365. " 2467 " − 740 " 2467.
366. A number is divisible by 2473 if − 1731 times the last digit of the number added to the rest of the number is divisible by 2473.
367. " 2477 " − 743 " 2477.
368. " 2503 " − 1752 " 2503.
369. " 2521 " − 252 " 2521.
370. " 2531 " − 253 " 2531.
371. A number is divisible by 2539 if − 2285 times the last digit of the number added to the rest of the number is divisible by 2539.
372. " 2543 " − 1780 " 2543.
373. " 2549 " − 2294 " 2549.
374. " 2551 " − 255 " 2551.
375. " 2557 " − 767 " 2557.
376. A number is divisible by 2579 if − 2321 times the last digit of the number added to the rest of the number is divisible by 2579.
377. " 2591 " − 259 " 2591.
378. " 2593 " − 1815 " 2593.
379. " 2609 " − 2348 " 2609.
380. " 2617 " − 785 " 2617.
381. A number is divisible by 2621 if − 262 times the last digit of the number added to the rest of the number is divisible by 2621.
382. " 2633 " − 1843 " 2633.
383. " 2647 " − 794 " 2647.
384. " 2657 " − 797 " 2657.
385. " 2659 " − 2393 " 2659.
386. A number is divisible by 2663 if − 1864 times the last digit of the number added to the rest of the number is divisible by 2663.
387. " 2671 " − 267 " 2671.
388. " 2677 " − 803 " 2677.
389. " 2683 " − 1878 " 2683.
390. " 2687 " − 806 " 2687.
391. A number is divisible by 2689 if − 2420 times the last digit of the number added to the rest of the number is divisible by 2689.
392. " 2693 " − 1885 " 2693.
393. " 2699 " − 2429 " 2699.
394. " 2707 " − 812 " 2707.
395. " 2711 " − 271 " 2711.
396. A number is divisible by 2713 if − 1899 times the last digit of the number added to the rest of the number is divisible by 2713.
397. " 2719 " − 2447 " 2719.
398. " 2729 " − 2456 " 2729.
399. " 2731 " − 273 " 2731.
400. " 2741 " − 274 " 2741.
401. A number is divisible by 2749 if − 2474 times the last digit of the number added to the rest of the number is divisible by 2749.
402. " 2753 " − 1927 " 2753.
403. " 2767 " − 830 " 2767.
404. " 2777 " − 833 " 2777.
405. " 2789 " − 2510 " 2789.
406. A number is divisible by 2791 if − 279 times the last digit of the number added to the rest of the number is divisible by 2791.
407. " 2797 " − 839 " 2797.
408. " 2801 " − 280 " 2801.
409. " 2803 " − 1962 " 2803.
410. " 2819 " − 2537 " 2819.
411. A number is divisible by 2833 if − 1983 times the last digit of the number added to the rest of the number is divisible by 2833.
412. " 2837 " − 851 " 2837.
413. " 2843 " − 1990 " 2843.
414. " 2851 " − 285 " 2851.
415. " 2857 " − 857 " 2857.
416. A number is divisible by 2861 if − 286 times the last digit of the number added to the rest of the number is divisible by 2861.
417. " 2879 " − 2591 " 2879.
418. " 2887 " − 866 " 2887.
419. " 2897 " − 869 " 2897.
420. " 2903 " − 2032 " 2903.

TABLE 1. SIMPLE DIVISIBILITY RULES FOR THE 1$^{st}$ 1000 PRIME NUMBERS—*Continued*

421. A number is divisible by 2909 if − 2618 times the last digit of the number added to the rest of the number is divisible by 2909.
422. " 2917 " − 875 " 2917.
423. " 2927 " − 878 " 2927.
424. " 2939 " − 2645 " 2939.
425. " 2953 " − 2067 " 2953.
426. A number is divisible by 2957 if − 887 times the last digit of the number added to the rest of the number is divisible by 2957.
427. " 2963 " − 2074 " 2963.
428. " 2969 " − 2672 " 2969.
429. " 2971 " − 297 " 2971.
430. " 2999 " − 2699 " 2999.
431. A number is divisible by 3001 if − 300 times the last digit of the number added to the rest of the number is divisible by 3001.
432. " 3011 " − 301 " 3011.
433. " 3019 " − 2717 " 3019.
434. " 3023 " − 2116 " 3023.
435. " 3037 " − 911 " 3037.
436. A number is divisible by 3041 if − 304 times the last digit of the number added to the rest of the number is divisible by 3041.
437. " 3049 " − 2744 " 3049.
438. " 3061 " − 306 " 3061.
439. " 3067 " − 920 " 3067.
440. " 3079 " − 2771 " 3079.
441. A number is divisible by 3083 if − 2158 times the last digit of the number added to the rest of the number is divisible by 3083.
442. " 3089 " − 2780 " 3089.
443. " 3109 " − 2798 " 3109.
444. " 3119 " − 2807 " 3119.
445. " 3121 " − 312 " 3121.
446. A number is divisible by 3137 if − 941 times the last digit of the number added to the rest of the number is divisible by 3137.
447. " 3163 " − 2214 " 3163.
448. " 3167 " − 950 " 3167.
449. " 3169 " − 2852 " 3169.
450. " 3181 " − 318 " 3181.
451. A number is divisible by 3187 if − 956 times the last digit of the number added to the rest of the number is divisible by 3187.
452. " 3191 " − 319 " 3191.
453. " 3203 " − 2242 " 3203.
454. " 3209 " − 2888 " 3209.
455. " 3217 " − 965 " 3217.
456. A number is divisible by 3221 if − 322 times the last digit of the number added to the rest of the number is divisible by 3221.
457. " 3229 " − 2906 " 3229.
458. " 3251 " − 325 " 3251.
459. " 3253 " − 2277 " 3253.
460. " 3257 " − 977 " 3257.
461. A number is divisible by 3259 if − 2933 times the last digit of the number added to the rest of the number is divisible by 3259.
462. " 3271 " − 327 " 3271.
463. " 3299 " − 2969 " 3299.
464. " 3301 " − 330 " 3301.
465. " 3307 " − 992 " 3307.
466. A number is divisible by 3313 if − 2319 times the last digit of the number added to the rest of the number is divisible by 3313.
467. " 3319 " − 2987 " 3319.
468. " 3323 " − 2326 " 3323.
469. " 3329 " − 2996 " 3329.
470. " 3331 " − 333 " 3331.
471. A number is divisible by 3343 if − 2340 times the last digit of the number added to the rest of the number is divisible by 3343.
472. " 3347 " − 1004 " 3347.
473. " 3359 " − 3023 " 3359.
474. " 3361 " − 336 " 3361.
475. " 3371 " − 337 " 3371.
476. A number is divisible by 3373 if − 2361 times the last digit of the number added to the rest of the number is divisible by 3373.
477. " 3389 " − 3050 " 3389.
478. " 3391 " − 339 " 3391.
479. " 3407 " − 1022 " 3407.
480. " 3413 " − 2389 " 3413.
481. A number is divisible by 3433 if − 2403 times the last digit of the number added to the rest of the number is divisible by 3433.
482. " 3449 " − 3104 " 3449.
483. " 3457 " − 1037 " 3457.
484. " 3461 " − 346 " 3461.
485. " 3463 " − 2424 " 3463.
486. A number is divisible by 3467 if − 1040 times the last digit of the number added to the rest of the number is divisible by 3467.
487. " 3469 " − 3122 " 3469.
488. " 3491 " − 349 " 3491.
489. " 3499 " − 3149 " 3499.
490. " 3511 " − 351 " 3511.
491. A number is divisible by 3517 if − 1055 times the last digit of the number added to the rest of the number is divisible by 3517.
492. " 3527 " − 1058 " 3527.
493. " 3529 " − 3176 " 3529.
494. " 3533 " − 2473 " 3533.
495. " 3539 " − 3185 " 3539.

TABLE 1. SIMPLE DIVISIBILITY RULES FOR THE 1$^{st}$ 1000 PRIME NUMBERS—*Continued*

496. A number is divisible by 3541 if − 354 times the last digit of the number added to the rest of the number is divisible by 3541.
497. " 3547 " − 1064 " 3547.
498. " 3557 " − 1067 " 3557.
499. " 3559 " − 3203 " 3559.
500. " 3571 " − 357 " 3571.
501. A number is divisible by 3581 if − 358 times the last digit of the number added to the rest of the number is divisible by 3581.
502. " 3583 " − 2508 " 3583.
503. " 3593 " − 2515 " 3593.
504. " 3607 " − 1082 " 3607.
505. " 3613 " − 2529 " 3613.
506. A number is divisible by 3617 if − 1085 times the last digit of the number added to the rest of the number is divisible by 3617.
507. " 3623 " − 2536 " 3623.
508. " 3631 " − 363 " 3631.
509. " 3637 " − 1091 " 3637.
510. " 3643 " − 2550 " 3643.
511. A number is divisible by 3659 if − 3293 times the last digit of the number added to the rest of the number is divisible by 3659.
512. " 3671 " − 367 " 3671.
513. " 3673 " − 2571 " 3673.
514. " 3677 " − 1103 " 3677.
515. " 3691 " − 369 " 3691.
516. A number is divisible by 3697 if − 1109 times the last digit of the number added to the rest of the number is divisible by 3697.
517. " 3701 " − 370 " 3701.
518. " 3709 " − 3338 " 3709.
519. " 3719 " − 3347 " 3719.
520. " 3727 " − 1118 " 3727.
521. A number is divisible by 3733 if − 2613 times the last digit of the number added to the rest of the number is divisible by 3733.
522. " 3739 " − 3365 " 3739.
523. " 3761 " − 376 " 3761.
524. " 3767 " − 1130 " 3767.
525. " 3769 " − 3392 " 3769.
526. A number is divisible by 3779 if − 3401 times the last digit of the number added to the rest of the number is divisible by 3779.
527. " 3793 " − 2655 " 3793.
528. " 3797 " − 1139 " 3797.
529. " 3803 " − 2662 " 3803.
530. " 3821 " − 382 " 3821.
531. A number is divisible by 3823 if − 2676 times the last digit of the number added to the rest of the number is divisible by 3823.
532. " 3833 " − 2683 " 3833.
533. " 3847 " − 1154 " 3847.
534. " 3851 " − 385 " 3851.
535. " 3853 " − 2697 " 3853.
536. A number is divisible by 3863 if − 2704 times the last digit of the number added to the rest of the number is divisible by 3863.
537. " 3877 " − 1163 " 3877.
538. " 3881 " − 388 " 3881.
539. " 3889 " − 3500 " 3889.
540. " 3907 " − 1172 " 3907.
541. A number is divisible by 3911 if − 391 times the last digit of the number added to the rest of the number is divisible by 3911.
542. " 3917 " − 1175 " 3917.
543. " 3919 " − 3527 " 3919.
544. " 3923 " − 2746 " 3923.
545. " 3929 " − 3536 " 3929.
546. A number is divisible by 3931 if − 393 times the last digit of the number added to the rest of the number is divisible by 3931.
547. " 3943 " − 2760 " 3943.
548. " 3947 " − 1184 " 3947.
549. " 3967 " − 1190 " 3967.
550. " 3989 " − 3590 " 3989.
551. A number is divisible by 4001 if − 400 times the last digit of the number added to the rest of the number is divisible by 4001.
552. " 4003 " − 2802 " 4003.
553. " 4007 " − 1202 " 4007.
554. " 4013 " − 2809 " 4013.
555. " 4019 " − 3617 " 4019.
556. A number is divisible by 4021 if − 402 times the last digit of the number added to the rest of the number is divisible by 4021.
557. " 4027 " − 1208 " 4027.
558. " 4049 " − 3644 " 4049.
559. " 4051 " − 405 " 4051.
560. " 4057 " − 1217 " 4057.
561. A number is divisible by 4073 if − 2851 times the last digit of the number added to the rest of the number is divisible by 4073.
562. " 4079 " − 3671 " 4079.
563. " 4091 " − 409 " 4091.
564. " 4093 " − 2865 " 4093.
565. " 4099 " − 3689 " 4099.
566. A number is divisible by 4111 if − 411 times the last digit of the number added to the rest of the number is divisible by 4111.
567. " 4127 " − 1238 " 4127.
568. " 4129 " − 3716 " 4129.
569. " 4133 " − 2893 " 4133.
570. " 4139 " − 3725 " 4139.

TABLE 1. SIMPLE DIVISIBILITY RULES FOR THE $1^{st}$ 1000 PRIME NUMBERS—*Continued*

571. A number is divisible by 4153 if − 2907 times the last digit of the number added to the rest of the number is divisible by 4153.
572. " 4157 " − 1247 " 4157.
573. " 4159 " − 3743 " 4159.
574. " 4177 " − 1253 " 4177.
575. " 4201 " − 420 " 4201.
576. A number is divisible by 4211 if − 421 times the last digit of the number added to the rest of the number is divisible by 4211.
577. " 4217 " − 1265 " 4217.
578. " 4219 " − 3797 " 4219.
579. " 4229 " − 3806 " 4229.
580. " 4231 " − 423 " 4231.
581. A number is divisible by 4241 if − 424 times the last digit of the number added to the rest of the number is divisible by 4241.
582. " 4243 " − 2970 " 4243.
583. " 4253 " − 2977 " 4253.
584. " 4259 " − 3833 " 4259.
585. " 4261 " − 426 " 4261.
586. A number is divisible by 4271 if − 427 times the last digit of the number added to the rest of the number is divisible by 4271.
587. " 4273 " − 2991 " 4273.
588. " 4283 " − 2998 " 4283.
589. " 4289 " − 3860 " 4289.
590. " 4297 " − 1289 " 4297.
591. A number is divisible by 4327 if − 1298 times the last digit of the number added to the rest of the number is divisible by 4327.
592. " 4337 " − 1301 " 4337.
593. " 4339 " − 3905 " 4339.
594. " 4349 " − 3914 " 4349.
595. " 4357 " − 1307 " 4357.
596. A number is divisible by 4363 if − 3054 times the last digit of the number added to the rest of the number is divisible by 4363.
597. " 4373 " − 3061 " 4373.
598. " 4391 " − 439 " 4391.
599. " 4397 " − 1319 " 4397.
600. " 4409 " − 3968 " 4409.
601. A number is divisible by 4421 if − 442 times the last digit of the number added to the rest of the number is divisible by 4421.
602. " 4423 " − 3096 " 4423.
603. " 4441 " − 444 " 4441.
604. " 4447 " − 1334 " 4447.
605. " 4451 " − 445 " 4451.
606. A number is divisible by 4457 if − 1337 times the last digit of the number added to the rest of the number is divisible by 4457.
607. " 4463 " − 3124 " 4463.
608. " 4481 " − 448 " 4481.
609. " 4483 " − 3138 " 4483.
610. " 4493 " − 3145 " 4493.
611. A number is divisible by 4507 if − 1352 times the last digit of the number added to the rest of the number is divisible by 4507.
612. " 4513 " − 3159 " 4513.
613. " 4517 " − 1355 " 4517.
614. " 4519 " − 4067 " 4519.
615. " 4523 " − 3166 " 4523.
616. A number is divisible by 4547 if − 1364 times the last digit of the number added to the rest of the number is divisible by 4547.
617. " 4549 " − 4094 " 4549.
618. " 4561 " − 456 " 4561.
619. " 4567 " − 1370 " 4567.
620. " 4583 " − 3208 " 4583.
621. A number is divisible by 4591 if − 459 times the last digit of the number added to the rest of the number is divisible by 4591.
622. " 4597 " − 1379 " 4597.
623. " 4603 " − 3222 " 4603.
624. " 4621 " − 462 " 4621.
625. " 4637 " − 1391 " 4637.
626. A number is divisible by 4639 if − 4175 times the last digit of the number added to the rest of the number is divisible by 4639.
627. " 4643 " − 3250 " 4643.
628. " 4649 " − 4184 " 4649.
629. " 4651 " − 465 " 4651.
630. " 4657 " − 1397 " 4657.
631. A number is divisible by 4663 if − 3264 times the last digit of the number added to the rest of the number is divisible by 4663.
632. " 4673 " − 3271 " 4673.
633. " 4679 " − 4211 " 4679.
634. " 4691 " − 469 " 4691.
635. " 4703 " − 3292 " 4703.
636. A number is divisible by 4721 if − 472 times the last digit of the number added to the rest of the number is divisible by 4721.
637. " 4723 " − 3306 " 4723.
638. " 4729 " − 4256 " 4729.
639. " 4733 " − 3313 " 4733.
640. " 4751 " − 475 " 4751.
641. A number is divisible by 4759 if − 4283 times the last digit of the number added to the rest of the number is divisible by 4759.
642. " 4783 " − 3348 " 4783.
643. " 4787 " − 1436 " 4787.
644. " 4789 " − 4310 " 4789.
645. " 4793 " − 3355 " 4793.

TABLE 1. SIMPLE DIVISIBILITY RULES FOR THE 1$^{st}$ 1000 PRIME NUMBERS—*Continued*

| | | | | | |
|---|---|---|---|---|---|
| 646. | A number is divisible by | 4799 if | − 4319 times the last digit of the number added to the rest of the number is divisible by | 4799. |
| 647. | ” | 4801 ” | − 480 | ” | 4801. |
| 648. | ” | 4813 ” | − 3369 | ” | 4813. |
| 649. | ” | 4817 ” | − 1445 | ” | 4817. |
| 650. | ” | 4831 ” | − 483 | ” | 4831. |
| 651. | A number is divisible by | 4861 if | − 486 times the last digit of the number added to the rest of the number is divisible by | 4861. |
| 652. | ” | 4871 ” | − 487 | ” | 4871. |
| 653. | ” | 4877 ” | − 1463 | ” | 4877. |
| 654. | ” | 4889 ” | − 4400 | ” | 4889. |
| 655. | ” | 4903 ” | − 3432 | ” | 4903. |
| 656. | A number is divisible by | 4909 if | − 4418 times the last digit of the number added to the rest of the number is divisible by | 4909. |
| 657. | ” | 4919 ” | − 4427 | ” | 4919. |
| 658. | ” | 4931 ” | − 493 | ” | 4931. |
| 659. | ” | 4933 ” | − 3453 | ” | 4933. |
| 660. | ” | 4937 ” | − 1481 | ” | 4937. |
| 661. | A number is divisible by | 4943 if | − 3460 times the last digit of the number added to the rest of the number is divisible by | 4943. |
| 662. | ” | 4951 ” | − 495 | ” | 4951. |
| 663. | ” | 4957 ” | − 1487 | ” | 4957. |
| 664. | ” | 4967 ” | − 1490 | ” | 4967. |
| 665. | ” | 4969 ” | − 4472 | ” | 4969. |
| 666. | A number is divisible by | 4973 if | − 3481 times the last digit of the number added to the rest of the number is divisible by | 4973. |
| 667. | ” | 4987 ” | − 1496 | ” | 4987. |
| 668. | ” | 4993 ” | − 3495 | ” | 4993. |
| 669. | ” | 4999 ” | − 4499 | ” | 4999. |
| 670. | ” | 5003 ” | − 3502 | ” | 5003. |
| 671. | A number is divisible by | 5009 if | − 4508 times the last digit of the number added to the rest of the number is divisible by | 5009. |
| 672. | ” | 5011 ” | − 501 | ” | 5011. |
| 673. | ” | 5021 ” | − 502 | ” | 5021. |
| 674. | ” | 5023 ” | − 3516 | ” | 5023. |
| 675. | ” | 5039 ” | − 4535 | ” | 5039. |
| 676. | A number is divisible by | 5051 if | − 505 times the last digit of the number added to the rest of the number is divisible by | 5051. |
| 677. | ” | 5059 ” | − 4553 | ” | 5059. |
| 678. | ” | 5077 ” | − 1523 | ” | 5077. |
| 679. | ” | 5081 ” | − 508 | ” | 5081. |
| 680. | ” | 5087 ” | − 1526 | ” | 5087. |
| 681. | A number is divisible by | 5099 if | − 4589 times the last digit of the number added to the rest of the number is divisible by | 5099. |
| 682. | ” | 5101 ” | − 510 | ” | 5101. |
| 683. | ” | 5107 ” | − 1532 | ” | 5107. |
| 684. | ” | 5113 ” | − 3579 | ” | 5113. |
| 685. | ” | 5119 ” | − 4607 | ” | 5119. |
| 686. | A number is divisible by | 5147 if | − 1544 times the last digit of the number added to the rest of the number is divisible by | 5147. |
| 687. | ” | 5153 ” | − 3607 | ” | 5153. |
| 688. | ” | 5167 ” | − 1550 | ” | 5167. |
| 689. | ” | 5171 ” | − 517 | ” | 5171. |
| 690. | ” | 5179 ” | − 4661 | ” | 5179. |
| 691. | A number is divisible by | 5189 if | − 4670 times the last digit of the number added to the rest of the number is divisible by | 5189. |
| 692. | ” | 5197 ” | − 1559 | ” | 5197. |
| 693. | ” | 5209 ” | − 4688 | ” | 5209. |
| 694. | ” | 5227 ” | − 1568 | ” | 5227. |
| 695. | ” | 5231 ” | − 523 | ” | 5231. |
| 696. | A number is divisible by | 5233 if | − 3663 times the last digit of the number added to the rest of the number is divisible by | 5233. |
| 697. | ” | 5237 ” | − 1571 | ” | 5237. |
| 698. | ” | 5261 ” | − 526 | ” | 5261. |
| 699. | ” | 5273 ” | − 3691 | ” | 5273. |
| 700. | ” | 5279 ” | − 4751 | ” | 5279. |
| 701. | A number is divisible by | 5281 if | − 528 times the last digit of the number added to the rest of the number is divisible by | 5281. |
| 702. | ” | 5297 ” | − 1589 | ” | 5297. |
| 703. | ” | 5303 ” | − 3712 | ” | 5303. |
| 704. | ” | 5309 ” | − 4778 | ” | 5309. |
| 705. | ” | 5323 ” | − 3726 | ” | 5323. |
| 706. | A number is divisible by | 5333 if | − 3733 times the last digit of the number added to the rest of the number is divisible by | 5333. |
| 707. | ” | 5347 ” | − 1604 | ” | 5347. |
| 708. | ” | 5351 ” | − 535 | ” | 5351. |
| 709. | ” | 5381 ” | − 538 | ” | 5381. |
| 710. | ” | 5387 ” | − 1616 | ” | 5387. |
| 711. | A number is divisible by | 5393 if | − 3775 times the last digit of the number added to the rest of the number is divisible by | 5393. |
| 712. | ” | 5399 ” | − 4859 | ” | 5399. |
| 713. | ” | 5407 ” | − 1622 | ” | 5407. |
| 714. | ” | 5413 ” | − 3789 | ” | 5413. |
| 715. | ” | 5417 ” | − 1625 | ” | 5417. |
| 716. | A number is divisible by | 5419 if | − 4877 times the last digit of the number added to the rest of the number is divisible by | 5419. |
| 717. | ” | 5431 ” | − 543 | ” | 5431. |
| 718. | ” | 5437 ” | − 1631 | ” | 5437. |
| 719. | ” | 5441 ” | − 544 | ” | 5441. |
| 720. | ” | 5443 ” | − 3810 | ” | 5443. |

TABLE 1. SIMPLE DIVISIBILITY RULES FOR THE 1$^{st}$ 1000 PRIME NUMBERS—*Continued*

721. A number is divisible by 5449 if − 4904 times the last digit of the number added to the rest of the number is divisible by 5449.
722. " 5471 " − 547 " 5471.
723. " 5477 " − 1643 " 5477.
724. " 5479 " − 4931 " 5479.
725. " 5483 " − 3838 " 5483.
726. A number is divisible by 5501 if − 550 times the last digit of the number added to the rest of the number is divisible by 5501.
727. " 5503 " − 3852 " 5503.
728. " 5507 " − 1652 " 5507.
729. " 5519 " − 4967 " 5519.
730. " 5521 " − 552 " 5521.
731. A number is divisible by 5527 if − 1658 times the last digit of the number added to the rest of the number is divisible by 5527.
732. " 5531 " − 553 " 5531.
733. " 5557 " − 1667 " 5557.
734. " 5563 " − 3894 " 5563.
735. " 5569 " − 5012 " 5569.
736. A number is divisible by 5573 if − 3901 times the last digit of the number added to the rest of the number is divisible by 5573.
737. " 5581 " − 558 " 5581.
738. " 5591 " − 559 " 5591.
739. " 5623 " − 3936 " 5623.
740. " 5639 " − 5075 " 5639.
741. A number is divisible by 5641 if − 564 times the last digit of the number added to the rest of the number is divisible by 5641.
742. " 5647 " − 1694 " 5647.
743. " 5651 " − 565 " 5651.
744. " 5653 " − 3957 " 5653.
745. " 5657 " − 1697 " 5657.
746. A number is divisible by 5659 if − 5093 times the last digit of the number added to the rest of the number is divisible by 5659.
747. " 5669 " − 5102 " 5669.
748. " 5683 " − 3978 " 5683.
749. " 5689 " − 5120 " 5689.
750. " 5693 " − 3985 " 5693.
751. A number is divisible by 5701 if − 570 times the last digit of the number added to the rest of the number is divisible by 5701.
752. " 5711 " − 571 " 5711.
753. " 5717 " − 1715 " 5717.
754. " 5737 " − 1721 " 5737.
755. " 5741 " − 574 " 5741.
756. A number is divisible by 5743 if − 4020 times the last digit of the number added to the rest of the number is divisible by 5743.
757. " 5749 " − 5174 " 5749.
758. " 5779 " − 5201 " 5779.
759. " 5783 " − 4048 " 5783.
760. " 5791 " − 579 " 5791.
761. A number is divisible by 5801 if − 580 times the last digit of the number added to the rest of the number is divisible by 5801.
762. " 5807 " − 1742 " 5807.
763. " 5813 " − 4069 " 5813.
764. " 5821 " − 582 " 5821.
765. " 5827 " − 1748 " 5827.
766. A number is divisible by 5839 if − 5255 times the last digit of the number added to the rest of the number is divisible by 5839.
767. " 5843 " − 4090 " 5843.
768. " 5849 " − 5264 " 5849.
769. " 5851 " − 585 " 5851.
770. " 5857 " − 1757 " 5857.
771. A number is divisible by 5861 if − 586 times the last digit of the number added to the rest of the number is divisible by 5861.
772. " 5867 " − 1760 " 5867.
773. " 5869 " − 5282 " 5869.
774. " 5879 " − 5291 " 5879.
775. " 5881 " − 588 " 5881.
776. A number is divisible by 5897 if − 1769 times the last digit of the number added to the rest of the number is divisible by 5897.
777. " 5903 " − 4132 " 5903.
778. " 5923 " − 4146 " 5923.
779. " 5927 " − 1778 " 5927.
780. " 5939 " − 5345 " 5939.
781. A number is divisible by 5953 if − 4167 times the last digit of the number added to the rest of the number is divisible by 5953.
782. " 5981 " − 598 " 5981.
783. " 5987 " − 1796 " 5987.
784. " 6007 " − 1802 " 6007.
785. " 6011 " − 601 " 6011.
786. A number is divisible by 6029 if − 5426 times the last digit of the number added to the rest of the number is divisible by 6029.
787. " 6037 " − 1811 " 6037.
788. " 6043 " − 4230 " 6043.
789. " 6047 " − 1814 " 6047.
790. " 6053 " − 4237 " 6053.
791. A number is divisible by 6067 if − 1820 times the last digit of the number added to the rest of the number is divisible by 6067.
792. " 6073 " − 4251 " 6073.
793. " 6079 " − 5471 " 6079.
794. " 6089 " − 5480 " 6089.
795. " 6091 " − 609 " 6091.

TABLE 1. SIMPLE DIVISIBILITY RULES FOR THE 1$^{st}$ 1000 PRIME NUMBERS—*Continued*

796. A number is divisible by 6101 if  − 610 times the last digit of the number added to the rest of the number is divisible by 6101.
797. " 6113 " − 4279 " 6113.
798. " 6121 " − 612 " 6121.
799. " 6131 " − 613 " 6131.
800. " 6133 " − 4293 " 6133.
801. A number is divisible by 6143 if − 4300 times the last digit of the number added to the rest of the number is divisible by 6143.
802. " 6151 " − 615 " 6151.
803. " 6163 " − 4314 " 6163.
804. " 6173 " − 4321 " 6173.
805. " 6197 " − 1859 " 6197.
806. A number is divisible by 6199 if − 5579 times the last digit of the number added to the rest of the number is divisible by 6199.
807. " 6203 " − 4342 " 6203.
808. " 6211 " − 621 " 6211.
809. " 6217 " − 1865 " 6217.
810. " 6221 " − 622 " 6221.
811. A number is divisible by 6229 if − 5606 times the last digit of the number added to the rest of the number is divisible by 6229.
812. " 6247 " − 1874 " 6247.
813. " 6257 " − 1877 " 6257.
814. " 6263 " − 4384 " 6263.
815. " 6269 " − 5642 " 6269.
816. A number is divisible by 6271 if − 627 times the last digit of the number added to the rest of the number is divisible by 6271.
817. " 6277 " − 1883 " 6277.
818. " 6287 " − 1886 " 6287.
819. " 6299 " − 5669 " 6299.
820. " 6301 " − 630 " 6301.
821. A number is divisible by 6311 if − 631 times the last digit of the number added to the rest of the number is divisible by 6311.
822. " 6317 " − 1895 " 6317.
823. " 6323 " − 4426 " 6323.
824. " 6329 " − 5696 " 6329.
825. " 6337 " − 1901 " 6337.
826. A number is divisible by 6343 if − 4440 times the last digit of the number added to the rest of the number is divisible by 6343.
827. " 6353 " − 4447 " 6353.
828. " 6359 " − 5723 " 6359.
829. " 6361 " − 636 " 6361.
830. " 6367 " − 1910 " 6367.
831. A number is divisible by 6373 if − 4461 times the last digit of the number added to the rest of the number is divisible by 6373.
832. " 6379 " − 5741 " 6379.
833. " 6389 " − 5750 " 6389.
834. " 6397 " − 1919 " 6397.
835. " 6421 " − 642 " 6421.
836. A number is divisible by 6427 if − 1928 times the last digit of the number added to the rest of the number is divisible by 6427.
837. " 6449 " − 5804 " 6449.
838. " 6451 " − 645 " 6451.
839. " 6469 " − 5822 " 6469.
840. " 6473 " − 4531 " 6473.
841. A number is divisible by 6481 if − 648 times the last digit of the number added to the rest of the number is divisible by 6481.
842. " 6491 " − 649 " 6491.
843. " 6521 " − 652 " 6521.
844. " 6529 " − 5876 " 6529.
845. " 6547 " − 1964 " 6547.
846. A number is divisible by 6551 if − 655 times the last digit of the number added to the rest of the number is divisible by 6551.
847. " 6553 " − 4587 " 6553.
848. " 6563 " − 4594 " 6563.
849. " 6569 " − 5912 " 6569.
850. " 6571 " − 657 " 6571.
851. A number is divisible by 6577 if − 1973 times the last digit of the number added to the rest of the number is divisible by 6577.
852. " 6581 " − 658 " 6581.
853. " 6599 " − 5939 " 6599.
854. " 6607 " − 1982 " 6607.
855. " 6619 " − 5957 " 6619.
856. A number is divisible by 6637 if − 1991 times the last digit of the number added to the rest of the number is divisible by 6637.
857. " 6653 " − 4657 " 6653.
858. " 6659 " − 5993 " 6659.
859. " 6661 " − 666 " 6661.
860. " 6673 " − 4671 " 6673.
861. A number is divisible by 6679 if − 6011 times the last digit of the number added to the rest of the number is divisible by 6679.
862. " 6689 " − 6020 " 6689.
863. " 6691 " − 669 " 6691.
864. " 6701 " − 670 " 6701.
865. " 6703 " − 4692 " 6703.
866. A number is divisible by 6709 if − 6038 times the last digit of the number added to the rest of the number is divisible by 6709.
867. " 6719 " − 6047 " 6719.
868. " 6733 " − 4713 " 6733.
869. " 6737 " − 2021 " 6737.
870. " 6761 " − 676 " 6761.

TABLE 1. SIMPLE DIVISIBILITY RULES FOR THE 1$^{st}$ 1000 PRIME NUMBERS—*Continued*

871. A number is divisible by 6763 if − 4734 times the last digit of the number added to the rest of the number is divisible by 6763.
872. ” 6779 ” − 6101 ” 6779.
873. ” 6781 ” − 678 ” 6781.
874. ” 6791 ” − 679 ” 6791.
875. ” 6793 ” − 4755 ” 6793.
876. A number is divisible by 6803 if − 4762 times the last digit of the number added to the rest of the number is divisible by 6803.
877. ” 6823 ” − 4776 ” 6823.
878. ” 6827 ” − 2048 ” 6827.
879. ” 6829 ” − 6146 ” 6829.
880. ” 6833 ” − 4783 ” 6833.
881. A number is divisible by 6841 if − 684 times the last digit of the number added to the rest of the number is divisible by 6841.
882. ” 6857 ” − 2057 ” 6857.
883. ” 6863 ” − 4804 ” 6863.
884. ” 6869 ” − 6182 ” 6869.
885. ” 6871 ” − 687 ” 6871.
886. A number is divisible by 6883 if − 4818 times the last digit of the number added to the rest of the number is divisible by 6883.
887. ” 6899 ” − 6209 ” 6899.
888. ” 6907 ” − 2072 ” 6907.
889. ” 6911 ” − 691 ” 6911.
890. ” 6917 ” − 2075 ” 6917.
891. A number is divisible by 6947 if − 2084 times the last digit of the number added to the rest of the number is divisible by 6947.
892. ” 6949 ” − 6254 ” 6949.
893. ” 6959 ” − 6263 ” 6959.
894. ” 6961 ” − 696 ” 6961.
895. ” 6967 ” − 2090 ” 6967.
896. A number is divisible by 6971 if − 697 times the last digit of the number added to the rest of the number is divisible by 6971.
897. ” 6977 ” − 2093 ” 6977.
898. ” 6983 ” − 4888 ” 6983.
899. ” 6991 ” − 699 ” 6991.
900. ” 6997 ” − 2099 ” 6997.
901. A number is divisible by 7001 if − 700 times the last digit of the number added to the rest of the number is divisible by 7001.
902. ” 7013 ” − 4909 ” 7013.
903. ” 7019 ” − 6317 ” 7019.
904. ” 7027 ” − 2108 ” 7027.
905. ” 7039 ” − 6335 ” 7039.
906. A number is divisible by 7043 if − 4930 times the last digit of the number added to the rest of the number is divisible by 7043.
907. ” 7057 ” − 2117 ” 7057.
908. ” 7069 ” − 6362 ” 7069.
909. ” 7079 ” − 6371 ” 7079.
910. ” 7103 ” − 4972 ” 7103.
911. A number is divisible by 7109 if − 6398 times the last digit of the number added to the rest of the number is divisible by 7109.
912. ” 7121 ” − 712 ” 7121.
913. ” 7127 ” − 2138 ” 7127.
914. ” 7129 ” − 6416 ” 7129.
915. ” 7151 ” − 715 ” 7151.
916. A number is divisible by 7159 if − 6443 times the last digit of the number added to the rest of the number is divisible by 7159.
917. ” 7177 ” − 2153 ” 7177.
918. ” 7187 ” − 2156 ” 7187.
919. ” 7193 ” − 5035 ” 7193.
920. ” 7207 ” − 2162 ” 7207.
921. A number is divisible by 7211 if − 721 times the last digit of the number added to the rest of the number is divisible by 7211.
922. ” 7213 ” − 5049 ” 7213.
923. ” 7219 ” − 6497 ” 7219.
924. ” 7229 ” − 6506 ” 7229.
925. ” 7237 ” − 2171 ” 7237.
926. A number is divisible by 7243 if − 5070 times the last digit of the number added to the rest of the number is divisible by 7243.
927. ” 7247 ” − 2174 ” 7247.
928. ” 7253 ” − 5077 ” 7253.
929. ” 7283 ” − 5098 ” 7283.
930. ” 7297 ” − 2189 ” 7297.
931. A number is divisible by 7307 if − 2192 times the last digit of the number added to the rest of the number is divisible by 7307.
932. ” 7309 ” − 6578 ” 7309.
933. ” 7321 ” − 732 ” 7321.
934. ” 7331 ” − 733 ” 7331.
935. ” 7333 ” − 5133 ” 7333.
936. A number is divisible by 7349 if − 6614 times the last digit of the number added to the rest of the number is divisible by 7349.
937. ” 7351 ” − 735 ” 7351.
938. ” 7369 ” − 6632 ” 7369.
939. ” 7393 ” − 5175 ” 7393.
940. ” 7411 ” − 741 ” 7411.
941. A number is divisible by 7417 if − 2225 times the last digit of the number added to the rest of the number is divisible by 7417.
942. ” 7433 ” − 5203 ” 7433.
943. ” 7451 ” − 745 ” 7451.
944. ” 7457 ” − 2237 ” 7457.
945. ” 7459 ” − 6713 ” 7459.

TABLE 1. SIMPLE DIVISIBILITY RULES FOR THE 1$^{st}$ 1000 PRIME NUMBERS—*Continued*

| # | Rule |
|---|---|
| 946. | A number is divisible by 7477 if $-2243$ times the last digit of the number added to the rest of the number is divisible by 7477. |
| 947. | " 7481 " $-748$ " 7481. |
| 948. | " 7487 " $-2246$ " 7487. |
| 949. | " 7489 " $-6740$ " 7489. |
| 950. | " 7499 " $-6749$ " 7499. |
| 951. | A number is divisible by 7507 if $-2252$ times the last digit of the number added to the rest of the number is divisible by 7507. |
| 952. | " 7517 " $-2255$ " 7517. |
| 953. | " 7523 " $-5266$ " 7523. |
| 954. | " 7529 " $-6776$ " 7529. |
| 955. | " 7537 " $-2261$ " 7537. |
| 956. | A number is divisible by 7541 if $-754$ times the last digit of the number added to the rest of the number is divisible by 7541. |
| 957. | " 7547 " $-2264$ " 7547. |
| 958. | " 7549 " $-6794$ " 7549. |
| 959. | " 7559 " $-6803$ " 7559. |
| 960. | " 7561 " $-756$ " 7561. |
| 961. | A number is divisible by 7573 if $-5301$ times the last digit of the number added to the rest of the number is divisible by 7573. |
| 962. | " 7577 " $-2273$ " 7577. |
| 963. | " 7583 " $-5308$ " 7583. |
| 964. | " 7589 " $-6830$ " 7589. |
| 965. | " 7591 " $-759$ " 7591. |
| 966. | A number is divisible by 7603 if $-5322$ times the last digit of the number added to the rest of the number is divisible by 7603. |
| 967. | " 7607 " $-2282$ " 7607. |
| 968. | " 7621 " $-762$ " 7621. |
| 969. | " 7639 " $-6875$ " 7639. |
| 970. | " 7643 " $-5350$ " 7643. |
| 971. | A number is divisible by 7649 if $-6884$ times the last digit of the number added to the rest of the number is divisible by 7649. |
| 972. | " 7669 " $-6902$ " 7669. |
| 973. | " 7673 " $-5371$ " 7673. |
| 974. | " 7681 " $-768$ " 7681. |
| 975. | " 7687 " $-2306$ " 7687. |
| 976. | A number is divisible by 7691 if $-769$ times the last digit of the number added to the rest of the number is divisible by 7691. |
| 977. | " 7699 " $-6929$ " 7699. |
| 978. | " 7703 " $-5392$ " 7703. |
| 979. | " 7717 " $-2315$ " 7717. |
| 980. | " 7723 " $-5406$ " 7723. |
| 981. | A number is divisible by 7727 if $-2318$ times the last digit of the number added to the rest of the number is divisible by 7727. |
| 982. | " 7741 " $-774$ " 7741. |
| 983. | " 7753 " $-5427$ " 7753. |
| 984. | " 7757 " $-2327$ " 7757. |
| 985. | " 7759 " $-6983$ " 7759. |
| 986. | A number is divisible by 7789 if $-7010$ times the last digit of the number added to the rest of the number is divisible by 7789. |
| 987. | " 7793 " $-5455$ " 7793. |
| 988. | " 7817 " $-2345$ " 7817. |
| 989. | " 7823 " $-5476$ " 7823. |
| 990. | " 7829 " $-7046$ " 7829. |
| 991. | A number is divisible by 7841 if $-784$ times the last digit of the number added to the rest of the number is divisible by 7841. |
| 992. | " 7853 " $-5497$ " 7853. |
| 993. | " 7867 " $-2360$ " 7867. |
| 994. | " 7873 " $-5511$ " 7873. |
| 995. | " 7877 " $-2363$ " 7877. |
| 996. | A number is divisible by 7879 if $-7091$ times the last digit of the number added to the rest of the number is divisible by 7879. |
| 997. | " 7883 " $-5518$ " 7883. |
| 998. | " 7901 " $-790$ " 7901. |
| 999. | " 7907 " $-2372$ " 7907. |
| 1000. | " 7919 " $-7127$ " 7919. |